\documentclass[10pt,twoside]{article}

\usepackage{amsmath}
\usepackage{amsfonts}
\usepackage{amssymb}
\usepackage{amscd}
\usepackage{amsthm}
\usepackage{amsbsy}
\usepackage{graphicx}
\usepackage{bm}

\usepackage{latexsym}

\def\@oddhead{\hfill \shorttitle \hfill \thepage}
\def\@evenhead{\thepage \hfill \shortauthor \hfill}
\def\@oddfoot{}
\def\@evenfoot{}
\numberwithin{equation}{section}

\def\whitebox{{\hbox{\hskip 1pt
 \vrule height 6pt depth 1.5pt
 \lower 1.5pt\vbox to 7.5pt{\hrule width
    3.2pt\vfill\hrule width 3.2pt}%
 \vrule height 6pt depth 1.5pt
 \hskip 1pt } }}
\def\qed{\ifhmode\allowbreak\else\nobreak\fi\hfill\quad\nobreak
     \whitebox\medbreak}
\newcommand{\pf}{\noindent{\bf Proof:}\ }

\newcommand{\ignore}[1]{}

\newcommand{\al}{\alpha}
\newcommand{\ep}{\epsilon}

\newcommand{\bZ}{{\mathbb Z}}

\newcommand{\bC}{{\mathbb C}}

\def\WW{\mathcal W}
\def\RR{\mathcal R}
\def\NN{\mathcal N}
\newcommand{\Zz}{{\mathbb Z}}

\newcommand{\gf}{{\mathbb F}}

\newcommand{\fq} {{\gf_{q}}}

\newcommand{\Tr}{\rm Tr}

\newtheorem{theorem}{Theorem}[section]

\newtheorem{lemma}[theorem]{Lemma}
\newtheorem{corollary}[theorem]{Corollary}

\newtheorem{remark}[theorem]{Remark}
\newtheorem{construction}[theorem]{Construction}

\date{Dedicated to Keqin Feng on his seventieth birthday}
\title{\ \\[0.4cm] \ \\ \bf  Partial difference sets from quadratic forms and $p$-ary weakly regular bent functions }
\author{Tao Feng\footnote{Department of Mathematical Sciences, University of Delaware, Newark, DE 19716, USA}\hspace{2mm}, Bin Wen\footnote{Department of Mathematics, Suzhou
University, Suzhou  215006, China and Department of Mathematics,
Changshu Institute of Technology, Changshu, Jiangsu 215500, China. }, Qing Xiang\footnote{Department of Mathematical Sciences, University of Delaware, Newark, DE 19716, USA. Research supported in part by NSF Grant DMS 1001557, and by the Overseas Cooperation Fund (grant 10928101) of China.}
and Jianxing Yin\footnote{Department of Mathematics, Suzhou
University, Suzhou  215006, China. Research supported in part by Natural Science Foundation of China (grant
10831002).}}
\begin{document}
%-------------------

\maketitle

%------------------------

\thispagestyle{empty}

%--------------------------------------
\begin{abstract}
\vskip 3mm\footnotesize{

\vskip 4.5mm
\noindent
We generalize  the construction of affine polar
graphs in two different ways to obtain new partial difference sets and
amorphic association schemes. The first generalization uses a combination 
of quadratic forms and uniform cyclotomy. In the second generalization we replace
the quadratic form in the affine polar graph construction by higher
degree homogeneous functions that are $p$-ary weakly regular bent. 
The negative Latin square type partial difference sets arising from the first generalization are new.

\vspace*{2mm}
\noindent{\bf 2000 Mathematics Subject Classification: 05E30, 05B10}

\vspace*{2mm}
\noindent{\bf Keywords and Phrases: Amorphic association scheme, association scheme, bent
function, difference set, $p$-ary bent function, partial difference
set, quadratic form, strongly regular graph, uniform cyclotomy.}}

\end{abstract}
\vspace*{-12.6cm}
%%%%%%%%%%%%%%%%%%%%%%%%%%%%%%%%%%%%%
\vspace*{2mm} \noindent\hspace*{82mm}
\begin{picture}(41,10)(0,0)\thicklines\setlength{\unitlength}{1mm}
\put(0,2){\line(1,0){41}} \put(0,16){\line(1,0){41}}%
\put(0,12.){\sl \copyright\hspace{1mm}Higher Education Press}
\put(0,7.8){\sl \hspace*{4.8mm}and International Press}%
\put(0,3.6){{\sl \hspace*{4.8mm}Beijing--Boston} }
\end{picture}

\vspace*{-17.6mm}\noindent{{\sl The title of\\
This book*****}\\SMM\,?, pp.\,1--?} \vskip8mm
%%%%%%%%%%%%%%%%%%
\vspace*{12.0cm}

%----------------------------------------------------------

%---------------------------------------------------------------

\section{Introduction}

Let $G$ be a finite (multiplicative) group of order $v$. A
$k$-element subset $D$ of $G$ is called a {\em $(v,k,\lambda)$
difference set\/} if the list of ``differences'' $xy^{-1}$, $x,y\in
D$, $x\neq y$, represents each nonidentity element in $G$ exactly
$\lambda$ times. Thus, $D$ is a $(v,k,\lambda)$ difference set in
$G$ if and only if it satisfies the following equation in the group
ring $\Zz[G]$:
\begin{equation}\label{DSeq}
DD^{(-1)}=(k-\lambda)1_G+\lambda G,
\end{equation}
where $D=\sum_{d\in D}d, D^{(-1)}=\sum_{d\in D}d^{-1}$,
$G=\sum_{g\in G}g$, and $1_G$ is the identity element of $G$. As an
example of difference sets, we mention the classical Paley
difference set (in the additive group of $\gf_q$) consisting of the
nonzero squares of $\gf_q$, where $q\equiv 3$ (mod 4). Difference
sets are the same objects as regular (i.e., sharply transitive)
symmetric designs. They are the subject of much study in the past 50
years. For a recent survey, see \cite{xiangsurvey}.

Again let $G$ be a finite (multiplicative) group of order $v$. A
$k$-element subset $D$ of $G$ is called a {\em $(v, k, \lambda,
\mu)$ partial difference set} (PDS, in short) provided that the list
of ``differences'' $xy^{-1}$, $x, y \in D$, $x \neq y$, contains
each nonidentity element of $D$ exactly $\lambda$ times and each
nonidentity element of $G \backslash D$ exactly $\mu$ times. Using
the group ring notation, we have that $D$ is a $(v, k, \lambda, \mu)$
partial difference set in $G$ if and only if
\begin{equation}\label{PDSeq}
DD^{(-1)}=\gamma 1_G+(\lambda-\mu)D+\mu G,
\end{equation}
where $\gamma=k-\mu$ if $1_G\not\in D$ and $\gamma=k-\lambda$ if
$1_G\in D$. A PDS with $\lambda=\mu$ is just a difference set. If
$D$ is a $(v, k, \lambda, \mu)$ PDS with $\lambda\neq \mu$, then
$D^{(-1)}=D$; in which case, equation~(\ref{PDSeq}) becomes
\begin{equation}\label{eq1.3}
D^2=\gamma 1_G+(\lambda-\mu)D+\mu G.
\end{equation}
A well-known example of PDS is the Paley PDS. Let $\gf_q$ be a
finite field of size $q$ with $q\equiv 1$ (mod 4). Then the set of
nonzero squares in $\gf_q$ forms a $(q,\frac {q-1}{2},
\frac{q-5}{4}, \frac{q-1}{4})$ PDS in the additive group of $\gf_q$,
which is called the {\em Paley PDS}.

Given a $(v, k, \lambda, \mu)$ partial difference set $D$ in $G$
with $1_G\not\in D$ and $D^{(-1)}=D$, one can construct a strongly
regular Cayley graph, ${\rm Cay}(G,D)$, whose vertex set is $G$, and
two vertices $x,y$ are adjacent if and only if $xy^{-1}\in D$. Such
a strongly regular graph ${\rm Cay}(G,D)$ has $G$ as a regular
automorphism group. For example, the strongly regular Cayley graph
constructed from the Paley PDS is the {\em Paley graph}. On the
other hand, if a strongly regular graph has a regular automorphism
group $G$, one can obtain a partial difference set in $G$. Therefore
partial difference sets are equivalent to strongly regular graphs
with a regular automorphism group. For a survey on partial
difference sets, we refere the reader to \cite{Ma94}. For connections among
partial difference sets, two-weight codes, projective
two-intersection sets, we refer the reader to \cite{CK86}. The following is a well-known construction of PDS. See for example
\cite{CK86}.

\begin{construction}\label{affpol}
Let $Q:\gf_q^{n}\rightarrow \gf_q$ be a nonsingular quadratic form,
where $n$ is even and $q$ is a power of an odd prime $p$, and let
$$D=\{x\in\gf_q^n\mid Q(x)\;{\mbox is\; a\; nonzero\; square}\}.$$
Then $D$ is a PDS in $(\gf_q^n,+)$. The corresponding strongly
regular graph ${\rm Cay}(\gf_q^n, D)$ is the so-called affine polar
graph.
\end{construction}

In this paper, we generalize Construction~\ref{affpol} in two
different directions. First we will replace the condition ``$Q(x)$ is a nonzero square"
in the above definition of $D$ by ``$Q(x)$ is a nonzero $e$th power
in $\gf_q$, where $e\geq 2$, $e|(q-1)$, and $p^j\equiv -1$ (mod $e$)
for some positive integer $j$". In doing so, we obtain new PDS
in $(\gf_q^n,+)$ and new amorphic association schemes on
$\gf_{q}^n$. We give the detailed statement of the our first 
generalization of Construction~\ref{affpol} below. 

\begin{theorem}\label{PU}
Let $p$ be a prime, $e\geq 2$, $q=p^{2j\gamma}$, where $\gamma\geq
1$, $e|(p^j+1)$ and $j$ is the smallest such positive integer. Let
$C_i$, $0\leq i\leq e-1$, be the cyclotomic classes of $\gf_q$ of
order $e$, $n=2m$ be an even positive integer, and $Q:
V=\gf_q^n\rightarrow \fq$ be a nonsingular quadratic form. Then each
of the sets $$D_{C_i}:=\{x \in V\mid Q(x)\in C_i\},\;0\leq i\leq
e-1,$$ is a PDS in $(V,+)$ with parameters $(N^2,(N-\epsilon)R,
\epsilon N+R^2-3\ep R, R^2-\ep R)$, where $N=q^{m}$, $R=fq^{m-1}$,
and $\ep=1$ or $-1$ according as $Q$ is hyperbolic or elliptic.
\end{theorem}

Second, we will replace the quadratic form in
Construction~\ref{affpol} by higher degree homogeneous functions
that are $p$-ary weakly regular bent (we will define the terms used
here in the sequel). In doing so we also generalize a recent
construction of PDS in \cite{TPF09} from the characteristic 3 case
to the case of arbitrary odd characteristic $p$. The detailed
statement of this generalization of Construction~\ref{affpol} is
postponed to Section 3 since it involves too many technical terms.

The description of our first generalization of
Construction~\ref{affpol} given above is quite straightforward
except that we did not explain why we impose the condition
$p^j\equiv -1$ (mod $e$). We now give the owed explanation.

Let $g$ be a fixed primitive element of $\fq$, and let $e\geq 2$ be
a divisor of $q-1$, and $f=(q-1)/e$. The $e$th {\em cyclotomic classes}
$C_0,C_1,\ldots, C_{e-1}$ are defined by
$$C_i=\{g^{i+ej}\mid 0\leq j\leq
f-1\},$$ where $0\leq i\leq e-1$. Let $\psi_{\al}, \al\in \fq$, be
the character of $(\fq,+)$ defined by
\begin{equation}\label{defchar}
\psi_{\al}(x)=\omega_p^{{\Tr}(\al
 x)},
\end{equation}
for all $x\in\gf_q$, where ${\Tr}: \gf_{p^n}\rightarrow \gf_p$ is the absolute trace
function. The {\em cyclotomic periods} $\eta_i$ of order $e$
are defined by
$$\eta_i=\sum_{z\in C_i}\psi_1(z),$$
where $i=0,1,\ldots ,e-1$.

It is clear that $\sum_{i=0}^{e-1}\eta_i=-1$. Therefore it is impossible to have $\eta_0=\eta_1=\cdots=\eta_{e-1}$ since, otherwise we have 
$\eta_i=-\frac{1}{e}$ for all $0\leq i\leq e-1$, implying that $\frac{1}{e},\; e\geq 2,$ is an algebraic integer, a contradiction. This fact motivates the following definition. 
We say that the cyclotomic periods $\eta_i$ are {\em uniform} if all but one are equal. It is well known \cite{BMW, bew} that
$\eta_i, 0\leq i\leq e-1$, are uniform if and only if $p^j\equiv -1$
(mod $e$) for some positive integer $j$ (here $p$ is the
characteristic of $\gf_q$). Therefore the essence in our first 
generalization of Construction~\ref{affpol} is to replace ``$Q(x)$
is a nonzero square" by ``$Q(x)$ is a nonzero $e$th power in
$\gf_q$, and $\eta_i$ are uniform".

Next we give the definitions of the terms used in our second generalization of Construction~\ref{affpol}. Let $p$ be a
prime, $n\geq 1$ be an integer, and $f: \gf_{p^{n}}\rightarrow
\gf_{p}$ be a function. The {\em Walsh coefficient} of $f$ at $b\in
\gf_{p^{n}}$ is defined by
$$
\WW_{f}(b)=\sum\limits_{x\in \gf_{p^{n}}}\omega_{p}^{f(x)+{\Tr}(bx)}
$$
where $\omega_p=e^{\frac{2\pi i}{p}}$ is a primitive complex
$p$th root of unity, and elements of $\gf_{p}$ are considered as
integers modulo $p$. The function $f$ is said to be {\em $p$-ary
bent} if $|\WW_{f}(b)|^{2}=p^{n}$ for all $b\in \gf_{p^n}$. A
$p$-ary bent function $f$ is said to be {\em regular} if for every
$b\in \gf_{p^{n}}$, $p^{-\frac{n}{2}}\WW_{f}(b)$ is equal to a
complex $p$-th root of unity. A $p$-ary bent function $f$ is said to
be {\em weakly regular} if there exists a complex number $u$ with
$|u| = 1$ such that $up^{-\frac{n}{2}}\WW_{f}(b) =
\omega_{p}^{f^{*}(b)}$ for all $b\in \gf_{p^{n}}$, where $f^*:
\gf_{p^{n}}\rightarrow \gf_{p}$ is a function. Clearly, a regular
$p$-ary bent function is weakly regular. Moreover, it is not
difficult to see that if $f$ is weakly regular bent, then so is
$f^{*}$; we call $f^*$ the {\em dual} of $f$. For example, let
$f:\gf_{p^n}\rightarrow \gf_p$ be a nonsingular quadratic form (here
we view $\gf_{p^n}$ as an $n$-dimensional vector space over
$\gf_p$), where $p$ is an odd prime. Then $f$ is a weakly regular
bent function (cf. \cite{HK06, FL07}).

Binary bent functions are usually called {\it Boolean bent
functions}, or simply {\it bent functions}. These functions were
first introduced by Rothaus \cite{Rothaus76} in 1976. Later Kumar,
Scholtz, and Welch \cite{KSW85} generalized the notion of a Boolean
bent function to that of a $p$-ary bent function.

Let $f:\gf_{2^{n}}\rightarrow \gf_{2}$ be a function. Then it is
well known \cite{dillon} that $f$ is bent if and only if
$D_1:=\{x\in\gf_{2^n}\mid f(x)=1\}$ is a difference set in
$(\gf_{2^{n}},+)$. Thus, given a bent function
$f:\gf_{2^{n}}\rightarrow \gf_{2}$, the inverse image of $1$
(respectively, $0$) is a difference set in $(\gf_{2^{n}},+)$. We
comment that a difference set $D$ in $(\gf_{2^{n}},+)$ is also a
partial difference set since $-D=D$.

Now let $p$ be an odd prime, and $f: \gf_{p^{n}}\rightarrow \gf_{p}$
be a $p$-ary bent function such that $f(-x)=f(x)$ for all $x\in
\gf_{p^n}$. For each $i\in\gf_p$, we define
$$D_i:=\{x\in \gf_{p^n}\mid f(x)=i\}.$$
It is then natural to ask whether $D_i$ is a PDS. It turns out that
the answer to this question is in general negative. But in a recent
paper \cite{TPF09}, Tan, Pott and Feng proved that under certain
conditions, the answer to the above question is indeed positive.
Specifically, let $f: \gf_{3^{2m}}\rightarrow \gf_{3}$ be a weakly
regular bent function such that $f(-x)=f(x)$ for all $x\in
\gf_{3^{2m}}$ and $f(0)=0$. Then it was shown in \cite{TPF09} that
$D_0\setminus\{0\}$, $D_1$ and $D_2$ are all partial difference sets
in $(\gf_{3^{2m}},+)$.

In this paper, we give a generalization of the aforementioned result
of Tan et al. To be precise, let $f: \gf_{p^{2m}}\rightarrow
\gf_{p}$ be a weakly regular bent function, where $p$ is an odd
prime. Define $D_{i}$ for each $i\in\gf_{p}$ as above. Let $\RR$
(respectively, $\NN$) denote the set of nonzero squares
(respectively, nonsquares) of $\gf_{p}$. If there exists a positive
integer $k$ satisfying $\gcd(k-1,p-1)=1$ such that $f(tx)=t^{k}f(x)$
for all $t\in\gf_{p}$ and all $x\in \gf_{p^{2m}}$, then
$D_0\setminus\{0\}$, $D_{\RR}:=\cup_{i\in \RR}D_{i}$ and
$D_{\NN}:=\cup_{i\in \NN}D_{i}$ are all partial difference sets in
$(\gf_{p^{2m}},+)$. This construction can also be viewed as a
generalization of Construction~\ref{affpol}.

The PDS constructed in this paper naturally lead to amorphic
association schemes, which we define below. Let $V$ be a finite set.
A $d$-class {\it symmetric association scheme} on $V$ is a
partition of $V\times V$ into sets $R_0$, $R_1, \ldots , R_d$
(called {\em relations, or associate classes}) such that

\begin{enumerate}
\item $R_0=\{(x,x)\mid x\in V\}$ (the diagonal relation);
\item $R_{\ell}$ is symmetric for $\ell=1,2,\ldots ,d$;
\item for all $i,j,k$ in $\{0,1,2,\ldots ,d\}$ there is an integer $p_{ij}^k$ such that, for all $(x,y)\in R_k$,
$$|\{z\in V \mid (x,z)\in R_i\; {\rm and}\; (z,y)\in R_j\}|=p_{ij}^k.$$
\end{enumerate}

Since each symmetric relation $R_\ell$, $1\leq \ell\leq d$,
corresponds to an undirected graph $G_\ell=(X, R_{\ell})$, $1\leq
\ell\leq d$, with vertex set $V$ and edge set $R_{\ell}$, we can
think of an association scheme $(V, \{R_\ell\}_{0\leq\ell\leq d})$
as an edge-decomposition of the complete graph on the vertex set $V$
into graphs $G_\ell$ on the same vertex set with the property that
for all $i,j,k$ in $\{1,2,\ldots ,d\}$ and for all $xy\in E(G_k)$,
$$|\{z\in V \mid xz\in E(G_i)\; {\rm and}\; zy\in E(G_j)\}|=p_{ij}^k,$$
where $E(G_k)$, $E(G_i)$, and $E(G_j)$ are the edge sets of $G_k$,
$G_i$ and $G_j$ respectively. The graphs $G_{\ell}$, $1\leq\ell\leq
d$, will be called {\it the graphs} of the association scheme $(V,
\{R_\ell\}_{0\leq \ell\leq d})$. A strongly regular graph (SRG) and
its complement form a symmetric 2-class association scheme. For more
background on association schemes and strongly regular graphs,
see~\cite{bcn}. Given an association scheme $(V, \{R_\ell\}_{0\leq
\ell\leq d})$, we can take the union of classes to form graphs with
larger edge sets (this process is called a {\em{fusion}}), but it is
not necessarily guaranteed that the fused collection of graphs will
form an association scheme on $V$. If an association scheme has the
property that any of its fusions is also an association scheme, then
we call the association scheme {\em{amorphic}}. A well-known and
important example of amorphic association schemes is given by the
cyclotomic association scheme on $\gf_q$ where the cyclotomic
periods are uniform \cite{BMW, BanMun}.

A $(v,k,\lambda,\mu)$ strongly regular graph is said to be of
{\em{Latin square type}} (respectively, {\it negative Latin square
type}) if $(v,k,\lambda,\mu) = (N^2, R(N-\epsilon), \epsilon N+R^2-3
\epsilon R, R^2 - \epsilon R)$ and $\epsilon = 1$ (respectively,
$\epsilon=-1$). If an association scheme is amorphic, then each of
its graphs is clearly strongly regular. Moreover, A. V. Ivanov
\cite{ivanov} showed that in an amorphic association scheme with at
least three classes, all graphs of the scheme are of Latin square
type, or all graphs are of negative Latin square type. The converse
of Ivanov's result is proved to be true in \cite{imy}. In fact even
more is true because Van Dam~\cite{vanDam} could prove the following
result.

\begin{theorem}
\label{vanDamtheorem} Let $V$ be a set of size $v$, let $\{G_1, G_2,
\ldots, G_d \}$ be an edge-decomposition of the complete graph on
$V$, where each $G_i$ is a strongly regular graph on $V$. If $G_i, 1
\leq i \leq d,$ are all of Latin square type or all of negative
Latin square type, then the decomposition is a $d$-class amorphic
association scheme on $V$.
\end{theorem}

Theorem~\ref{vanDamtheorem} will be used in Section 2 to show that
the PDS constructed there actually lead to amorphic
association schemes. For a recent survey of results on amorphic
association schemes, we refer the reader to \cite{vdMuzy}.

The rest of the paper is organized as follows. In Sections 2 and 3
we give the two generalizations of Construction~\ref{affpol}. In
each of those two sections, we show that the PDS obtained actually
lead to amorphic association schemes. As far as we know, the PDS
with negative Latin square type parameters constructed in Section 2 are new.

\section{Partial difference sets from quadratic forms and uniform cyclotomy}

We start with some background on characters of abelian groups. Let
$G$ be a finite abelian group. A (complex) character $\chi$ of $G$
is a homomorphism from $G$ to $\bC^{*}$, the multiplicative group of
$\bC$. A character $\chi$ of $G$ is called {\it principal} if
$\chi(g) = 1$ for all $g\in G$; otherwise it is called {\it
nonprincipal}. The set of all characters of $G$ forms a group (under
point-wise multiplication), which is isomorphic to $G$. Let $\chi$
be a character of $G$, and $A=\sum\limits_{g\in G}a_{g}g\in\bC[G]$.
We define
$$\chi(A)=\sum\limits_{g\in G}a_g\chi(g).$$
Starting with the important work of Turyn \cite{turyn}, character
sums have been a powerful tool in the study of difference sets of
all types.  The following lemma states how character sums can be
used to verify that a subset of an abelian group is a PDS.

\begin{lemma}\label{CV}
Let $G$ be a multiplicatively written abelian group of order $v$
and $D$ be a subset of $G$ such that $D=\{d^{-1}\mid d \in D\}$, and
$1_G\not\in D$. Let $k,\lambda,\mu$ be positive integers such that
$k^2=\mu v+(\lambda-\mu)k+(k-\mu)$. Then $D$ is a
$(v,k,\lambda,\mu)$ PDS in $G$ if and only if
\[\chi(D)=\begin{cases}k,\quad &\textup{if $\chi$ is
principal on
$G$},\\\frac{(\lambda-\mu)\pm\sqrt{(\lambda-\mu)^2+4(k-\mu)}}{2},\quad
&\textup{if $\chi$ is nonprincipal on $G$}. \end{cases} \]
\end{lemma}

Let $q$ be a power of a prime $p$, and $e\geq 2$ be a divisor of
$q-1$. Let $g$ be a fixed primitive element of $\gf_q$, $C_0=\langle g^e\rangle$, $C_i=g^iC_0$, $1\leq i\leq e-1$, be the cyclotomic classes of
$\gf_q$ of order $e$, and $\eta_i=\sum_{x\in C_i}\psi_1(x)$ be the
cyclotomic periods of order $e$, where $\psi_1$ is defined in
(\ref{defchar}). The cyclotomic periods $\eta_i$ are in general very
difficult to determine when $e$ is large. But in a special case, one
can compute $\eta_i$ explicitly.

\begin{theorem}\label{UC}
Let $p$ be a prime, $e\geq 2$, $q=p^{2j\gamma}$, where $\gamma\geq
1$, $e|(p^j+1)$ and $j$ is the smallest such positive integer. Then
the cyclotomic periods are given by

{\bf Case A.} If $\gamma$, $p$, $\frac{p^j+1}{e}$ are all odd, then
\[\eta_{e/2}=\sqrt{q}-\frac{\sqrt{q}+1}{e},\;  \eta_i=-\frac{1+\sqrt{q}}{e},\; {\mbox for \;all} \;i\neq \frac{e}{2}.\]

{\bf Case B.} In all the other cases,
\[\eta_0=-(-1)^{\gamma}\sqrt{q}+\frac{(-1)^{\gamma}\sqrt{q}-1}{e},\; \eta_i=\frac{(-1)^{\gamma}\sqrt{q}-1}{e},\; {\mbox for \;all} \; i\neq 0.\]
\end{theorem}

For a proof of Theorem~\ref{UC}, we refer the reader to \cite{BMW,
MY}.

Next we define what we mean by nonsingular quadratic forms. Let $V$ be an $n$-dimensional vector space over $\gf_q$.
A function $Q: V \rightarrow \gf_q$ is called a {\it quadratic form}
if
\begin{enumerate}
\item  $Q(\al v) = \al^2 Q(v)$ for all $\al \in
\gf_q$ and $v \in V$,
\item the function $B: V \times V \rightarrow \gf_q$ defined by
$B(v_1,v_2) = Q(v_1+v_2) - Q(v_1) - Q(v_2)$ is bilinear.
\end{enumerate}
We say that  $Q$ is {\em{nonsingular}} if the subspace $W$ of $V$ with the
property that $Q$ vanishes on $W$ and $B(w,v)=0$ for all $v\in V$
and $w\in W$ is the zero subspace. If the field $\gf_q$ has odd
characteristic, then $Q$ is nonsingular if and only if $B$ is
nondegenerate; but this may not be true when $\gf_q$ has
characteristic 2, because in that case $Q$ may not be zero on the
radical ${\rm Rad}(V)=\{w\in V\mid B(w,v)=0$ for all $v\in V\}$.
However, if $V$ is an even-dimensional vector space over an
even-characteristic field $\gf_q$, then $Q$ is nonsingular if and
only if $B$ is nondegenerate (cf. \cite[p.~14]{cameron}).

Now assume that $n=2m$ is an even positive integer, and $Q:
V=\gf_q^n\rightarrow \fq$ is a nonsingular quadratic form. Therefore
the polar form $B(x,y):=Q(x+y)-Q(x)-Q(y)$ of $Q$ is nondegenerate.
We write $\chi_b$, $b\in V$, for the additive character of $V$
defined by
$$\chi_b(x)=\psi_1(B(b,x)),\; x\in V.$$ Since $B$ is
nondegenerate, we see that $\{\chi_b\mid b\in V\}$ is the set of all
additive characters of $V$.

For each $u\in \fq$, we define $D_{u}=\{x \in V\mid Q(x)=u\}$, and $\mathcal{L}_{\alpha}=\sum_{u\in \fq}
D_{u}\psi_{\alpha}(u)$ for each $\al\in \fq$, where $\psi_{\al}$ is
defined in (\ref{defchar}).  (Here $\mathcal{L}_{\alpha}$ is view as an element of the group ring $\bC[V,+]$.) For a subset $X$ of $\fq$, we write
$$D_X=\sum_{x\in X}D_x, \;\mathcal{L}_X=\sum_{x\in
X}\mathcal{L}_x.$$ Our main result in this section is
Theorem~\ref{PU}, as stated in Section 1. Below we give the proof of
that theorem.

\vspace{0.1in}

\noindent{\bf Proof of Theorem 1.2.} We will compute the character
sums $\chi_b(D_{C_i})$, $b\in V$, explicitly, and then use
Lemma~\ref{CV} to finish the proof. To simplify notation, we write
$C$ for $C_0$.

First, we have \begin{align*}\mathcal{L}_C=&\sum_{\al\in C}\sum_{u\in \fq}D_u\psi_{\al}(u)\\
                     =&fD_0+\sum_{\al\in C}\sum_{i=0}^{e-1}\sum_{u\in C_i}D_u\psi_{\al}(u)
                     =fD_0+\sum_{i=0}^{e-1}D_{C_i}\psi_1(C_i)
                       \end{align*}
Corresponding to the two cases of Theorem~\ref{UC}, we have
\begin{equation}\label{expLC}
\mathcal{L}_C=\begin{cases}(f-\eta_0)D_0+D_{C_{e/2}}(\eta_{e/2}-\eta_0)+\eta_0
V,\quad &\textup{in Case
A},\\(f-\eta_1)D_0+D_{C}(\eta_{0}-\eta_1)+\eta_1 V,\quad &\textup{in
Case B}.\end{cases}
\end{equation}

Given any $b\in V$, we now compute $\chi_{b}(\mathcal{L}_C)$.
\begin{align*}\chi_b(\mathcal{L}_C)=&\sum_{\al\in C}\sum_{u\in \fq}\chi_b(D_u)\psi_{\al}(u)=\sum_{\al\in C}\sum_{x\in V}\chi_b(x)\psi_{\al}(Q(x))\\=&\sum_{\al\in C}\sum_{x\in V}\psi_1(\al Q(x)+B(b,x))
=\sum_{\al\in C}\sum_{x\in V}\psi_1(\al Q(x+\al^{-1}b)-\al^{-1} Q(b))\\
=&\sum_{\al\in C}\psi_1(-\al^{-1} Q(b))\sum_{x\in V}\psi_1(\al
Q(x)).
\end{align*}
For each $\al \in \gf_q^*$, $\al Q(x)$ is a quadratic form which is
nonsingular and of the same type as $Q(x)$. By \cite[Theorem
3.2]{leep}, we have for each $\al \in \gf_q^*$, $\sum_{x\in
V}\psi_1(\al Q(x))=\epsilon q^{m}$, where $\epsilon=1$ if $Q$ is
hyperbolic and $\epsilon=-1$ if $Q$ is elliptic. Therefore, in Case
A we have
\begin{equation*}\chi_b(\mathcal{L}_C)=\ep q^{m}\sum_{\al\in C}\psi_1(-\al^{-1} Q(b))=\begin{cases}\ep q^{m}f, \quad\textup{if $Q(b)=0$},\\ \ep q^{m}\eta_{e/2}, \quad\textup{if $-Q(b)\in C_{e/2}$},\\ \ep q^{m}\eta_0, \quad\textup{otherwise},\end{cases}\end{equation*}
and, in Case B we have $$\chi_b(\mathcal{L}_C)=\begin{cases}\ep
q^{m}f, \quad\textup{if $Q(b)=0$},\\ \ep q^{m}\eta_0,
\quad\textup{if $-Q(b)\in C$},\\ \ep q^{m}\eta_1,
\quad\textup{otherwise}.\end{cases}$$

Next we compute $\chi_b(D_0)$. We have

\begin{align*}q\chi_b(D_0)=&\sum_{x\in V}\sum_{u\in\fq}\chi_b(x)\psi_u(Q(x))=\sum_{x\in V}\sum_{u\in\fq}\psi_1(B(b,x)+uQ(x))\\
=&\sum_{x\in V}\psi_1(B(b,x))+\sum_{x\in
V}\sum_{u\in\gf_q^*}\psi_1(B(b,x)+uQ(x)).
\end{align*}

We now restrict our attention to the case where $b\neq 0$;  in that
case, we have $\sum_{x\in V}\psi_1(B(b,x))=0$. It follows that

\begin{align*}
q\chi_b(D_0)=&\sum_{x\in V}\sum_{u\in\gf_q^*}\psi_1(-u^{-1}Q(b)+uQ(x+u^{-1}b)))\\
=&\sum_{u\in\gf_q^*}\psi_1(-u^{-1}Q(b))\sum_{x\in V}\psi_1(uQ(x))
=\ep q^{m}\sum_{u\in\gf_q^*}\psi_1(-u^{-1}Q(b))\\
=&\begin{cases}\ep q^{m} (q-1),&\quad \textup{if $Q(b)=0$},\\
  -\ep q^{m},& \quad \textup{otherwise}.
 \end{cases}
\end{align*}

From (\ref{expLC}), in Case A, we have
$$\chi_b(D_{C_{e/2}})=\frac{\chi_b(\mathcal{L}_C)-(f-\eta_0)\chi_b(D_0)}{\eta_{e/2}-\eta_0}.$$
Substituting the values of $\chi_b(\mathcal{L}_C)$ and
$\chi_b(D_0)$, we obtain the following: for
each $b\in V$, $b\neq 0$,
\begin{equation*}\chi_b(D_{C_{e/2}})=\begin{cases}\ep q^{m-1}\frac{q\eta_{e/2}-\eta_0+f}{\eta_{e/2}-\eta_0}=\ep q^{m-1}(q-f),\quad &\textup{if $-Q(b)\in C_{e/2}$},\\ \ep q^{m-1}\frac{(q-1)\eta_{0}+f}{\eta_{e/2}-\eta_0}=-\ep q^{m-1}f,\quad &\textup{otherwise}.\end{cases}
\end{equation*}

Again, from (\ref{expLC}), in Case B, we have
$$\chi_b(D_C)=\frac{\chi_b(\mathcal{L}_C)-(f-\eta_1)\chi_b(D_0)}{\eta_0-\eta_1}.$$
Substituting the values of $\chi_b(\mathcal{L}_C)$ and
$\chi_b(D_0)$, we obtain the following: for
each $b\in V$, $b\neq 0$,
\begin{equation*}\chi_b(D_C)=\begin{cases}\ep
q^{m-1}\frac{q\eta_0-\eta_1+f}{\eta_0-\eta_1}=\ep q^{m-1}(q-f),\quad
&\textup{if $-Q(b)\in C$},\\ \ep
q^{m-1}\frac{(q-1)\eta_1+f}{\eta_0-\eta_1}=-\ep q^{m-1}f,\quad
&\textup{otherwise}.\end{cases}
\end{equation*}

The sizes of $D_{C_{e/2}}$ and $D_C$ can be computed as follows.
First we have $|D_0|=\chi_0(D_0)=q^{2m-1}+\ep q^{m-1} (q-1)$. By
applying the principal character $\chi_0$ to ${\mathcal
L}_C=(e-\eta_0)D_0+D_{C_{e/2}}(\eta_{e/2}-\eta_0)+\eta_0 V$, we
solve that $|D_{C_{e/2}}|=(q^{m}-\epsilon)fq^{m-1}$ in Case A.
Similarly, we have $|D_C|=(q^{m}-\epsilon)fq^{m-1}$ in Case B.

By Lemma \ref{CV} and the computed character values of $D_{C_{e/2}}$
and $D_C$, we conclude that $D_{C_{e/2}}$ is a PDS in $(V,+)$ in
Case A, $D_{C}$ is a PDS in Case B. Both of them have parameters
$(N^2,(N-\epsilon)R, \epsilon N+R^2-3\ep R, R^2-\ep R)$, where
$N=q^{m}$, $R=fq^{m-1}$.

Now to prove that each $D_{C_i}$ is a PDS in $(V,+)$, one only needs
to replace the quadratic form $Q$ by $\alpha Q$ for an appropriate
$\alpha\in \gf_q^*$ and note that $\alpha Q$ is also a nonsingular
quadratic form on $V$ and of the same type as $Q$.

The proof is now complete.\qed

\begin{remark}
(1). It is well known that $D_0\setminus\{0\}$ is a PDS with
parameters $(q^{2m}, (q^m-\epsilon)r, q^m+r^2-3r, r^2-r)$, where
$r=q^{m-1}+\epsilon$, $\epsilon=1$ or $-1$ according as $Q$ is
hyperbolic or elliptic. This PDS is referred to as Example {\bf RT2}
in \cite{CK86}.

(2). It is interesting to note that Theorem~\ref{PU} is valid when
$p=2$ while Construction~\ref{affpol} only works when $p$ is odd.

(3). When the quadratic form $Q$ in Theorem~\ref{PU} is of elliptic
type, the PDS $D_{C_i}$, $0\leq i\leq (e-1)$, have negative Latin
square type parameters. Negative Latin square type PDS are harder to
come by than Latin square type PDS. Besides Example {\bf RT2} and
the PDS arising from Construction~\ref{affpol}, there is one more
general class of negative Latin square type PDS in elementary
abelian $p$-groups coming from the ``difference of two quadrics"
construction in \cite{brouwer}. (There exist variations of Brouwer's
construction, see \cite{ham, davisX}.) The negative Latin square
type PDS arising from Theorem~\ref{PU} have very different
parameters from those in \cite{brouwer, ham, davisX} since there is
quite a bit of freedom in choosing the parameter $f=(q-1)/e$. As far
as we know, the negative Latin square type PDS from Theorem~\ref{PU}
are new.
\end{remark}

Next we show that the PDS obtained in Theorem~\ref{PU} also lead to
amorphic association schemes.

\begin{corollary}\label{amor}
Let $p$ be a prime, $e\geq 2$, $q=p^{2j\gamma}$, where $\gamma\geq
1$, $e|(p^j+1)$ and $j$ is the smallest such positive integer. Let
$C_i$, $0\leq i\leq e-1$, be the cyclotomic classes of $\gf_q$ of
order $e$, $n=2m$ be an even positive integer, and $Q:
V=\gf_q^n\rightarrow \fq$ be a nonsingular quadratic form. Then the
Cayley graphs ${\rm Cay}(G,D_{0}\setminus\{0\})$, ${\rm
Cay}(G,D_{C_i})$, $0\leq i\leq e-1$, where $G=(V,+)$, form a
$(e+1)$-class amorphic association scheme on $V$.
\end{corollary}

\begin{Proof}
The proof is immediate by combining Theorem~\ref{PU} and
Theorem~\ref{vanDamtheorem}.\qed
\end{Proof}

\begin{remark}
From Corollary~\ref{amor} one obtains more PDS by choosing an
arbitrary subset of $\{D_{0}\setminus\{0\}, D_{C_i}: 0\leq i\leq
e-1\}$, and taking union of the members of the subset.
\end{remark}

\section{Partial difference sets from weakly regular $p$-ary bent functions}

Throughout this section, $p$ always denotes an {\bf odd} prime,
$p^*=(-1)^{\frac{p-1}{2}}p$, and $\omega_p=e^{\frac{2\pi i}{p}}$.
Let $\RR$ (respectively, $\NN$) denote the set of nonzero squares
(respectively, nonsquares) of $\gf_{p}$. From the values of
quadratic Gauss sums \cite[p.~22]{bew}, we have
$$r_0:=\sum_{a\in \RR}\omega_p^a=\frac{\sqrt{p^*} -1}{2},$$
and
$$n_0:=\sum_{a\in \NN}\omega_p^a=\frac{-\sqrt{p^*} -1}{2}.$$

Let $f: \gf_{p^{n}}\rightarrow \gf_{p}$ be a weakly regular bent
function. From \cite{FL07,HK06}, we have
\begin{equation}\label{WRf}
{\mathcal{W}}_{f}(b)=u(\sqrt{p^{*}})^{n}\omega_{p}^{f^{*}(b)},
\end{equation}
for every $b\in\gf_{p^n}$, where $u=\pm1$ and $f^*$ is a function
from $\gf_{p^n}$ to $\gf_p$. Furthermore, we define
\begin{equation}\label{defD}
D_{i}=\{x\in \gf_{p^{n}}|\ f(x)=i\},
\end{equation}
for every $i\in\gf_{p}$, and
$${\mathcal{L}}_{t}=\sum\limits_{i=0}^{p-1}D_{i}\omega_{p}^{it},$$ for
every $t\in\gf_{p}$. (Here ${\mathcal{L}}_{t}$ is viewed as an
element of the group ring $\bC[(\gf_{p^n},+)]$). Furthermore, we
define
$$D_{\RR}=\cup_{i\in\RR}D_i,\;D_{\NN}=\cup_{i\in\NN}D_i,$$
and
$${\mathcal{L}}_{\RR}=\sum_{t\in\RR}{\mathcal{L}}_t, \;
{\mathcal{L}}_{\NN}=\sum_{t\in\NN}{\mathcal{L}}_t.$$

We will need a couple of lemmas from \cite{PTFL}.
\begin{lemma} {\em (\cite{PTFL})}\label{dualdeg}
Let $f:\gf_{p^{n}}\rightarrow \gf_{p}$ be a weakly regular bent
function. Assume that the Walsh coefficients of $f$ satisfy
(\ref{WRf}). If there exists a constant $k$ with $\gcd(k-1,p-1)=1$
such that $f(tx)=t^{k}f(x)$ for all $t\in\gf_{p}$ and
 all $x\in \gf_{p^{n}}$, then there exists a constant $\ell$ satisfying
$\gcd(\ell-1,p-1)=1$ and $f^{*}(tx)=t^{\ell}f^{*}(x)$ for all
$t\in\gf_{p}$ and all $x\in \gf_{p^{n}}$.
\end{lemma}

\begin{remark} The condition that there exists some $k$ with $\gcd(k-1,p-1)=1$
such that $f(tx)=t^{k}f(x)$ for all $t\in\gf_{p}$ and
 all $x\in \gf_{p^{n}}$ is not a severe one. Almost all known weakly regular bent functions satisfy this condition.
From the proof in \cite{PTFL}, we see that the constant $\ell$ can
be chosen in such a way that $\ell\equiv k(k-1)^{-1}$ (mod $p$).
\end{remark}

\begin{lemma} {\em (\cite{PTFL})}\label{ProdL}
Let $f:\gf_{p^{n}}\rightarrow \gf_{p}$ be a weakly regular bent
function. Assume that the Walsh coefficients of $f$ satisfy
(\ref{WRf}) and there exists a constant $k$ with $\gcd(k-1,p-1)=1$
such that $f(tx)=t^{k}f(x)$ for all $t\in\gf_{p}$ and all $x\in
\gf_{p^{n}}$. Let $\ell$ be given by Lemma~\ref{dualdeg}. Then
\begin{enumerate}
\item [(1)] for $s,t, s+t\in\gf_{p}^*$,
${\mathcal{L}}_{t}{\mathcal{L}}_{s}=u\bigg(\frac{tsv}{p}\bigg)^{n}(\sqrt{p^{*}})^{n}{\mathcal{L}}_{v},$
where $v=(s^{1-\ell}+t^{1-\ell})^{\frac{1}{1-\ell}}$;
\item [(2)] for $t\in\gf_{p}^{*}$,
${\mathcal{L}}_{t}{\mathcal{L}}_{-t}=p^{n}$;
\item [(3)] for $a\in\gf_{p}$,
$\sum\limits_{t=1}^{p-1}{\mathcal{L}}_{t}{\mathcal{L}}_{0}\omega_{p}^{-at}=(p|D_{a}|-p^{n})\gf_{p^{n}}$.
\end{enumerate}
\end{lemma}

\begin{remark}
The symbol $\bigg(\frac{\cdot}{p}\bigg)$ in part (1) of
Lemma~\ref{ProdL} is the Lengdre symbol. The equalities in all three
parts of the lemma should be viewed as equalities in the group ring
$\bC[(\gf_{p^n},+)]$. In particular, the right hand side of the
equality in part (2) of Lemma~\ref{ProdL} is really $p^n\cdot 0$,
where $0$ is the identity of $(\gf_{p^n},+)$.

\end{remark}

Now we are able to establish the main result of this section.

\begin{theorem}\label{SC1}
Let $f:\gf_{p^{n}}\rightarrow \gf_{p}$ be a weakly regular bent
function. Assume that the Walsh coefficients of $f$ satisfy
(\ref{WRf}) and there exists a constant $k$ with $\gcd(k-1,p-1)=1$
such that $f(tx)=t^{k}f(x)$ for all $t\in\gf_{p}$ and all $x\in
\gf_{p^{n}}$. Furthermore, assume that $n$ is even. Then
$D_{0}\setminus \{0\}$, $D_{\RR}$ and $D_{\NN}$ are all partial
difference sets in $(\gf_{p^{n}},+)$.
\end{theorem}

\pf By assumption we have $f(-x)=f(x)$ for all $x\in\gf_{p^n}$. It
follows that $0\in D_0$, $-D_0=D_0$, $-D_{\RR}=D_{\RR}$, and
$-D_{\NN}=D_{\NN}$. We will compute $D_0^2$, $D_{\RR}^2$ and
$D_{\NN}^2$ in the group ring $\bC[(\gf_{p^n},+)]$. To this end, we
first make some preparation.

By definition, we have
\begin{eqnarray*}
{\mathcal{L}}_{\RR}&=&\sum_{t\in\RR}{\mathcal{L}}_t=\sum_{t\in\RR}\sum_{i=0}^{p-1}D_i\omega_p^{it}\\
                   &=&\sum_{i=0}^{p-1}D_i\bigg(\sum_{t\in\RR}\omega_p^{it}\bigg)\\
                   &=&\frac{p-1}{2}D_0+r_0D_{\RR}+n_0D_{\NN}
\end{eqnarray*}
Similarly,
${\mathcal{L}}_{\NN}=\frac{p-1}{2}D_0+n_0D_{\RR}+r_0D_{\NN}$. Also
from definition, we have ${\mathcal L}_0=\gf_{p^n}$.

From the definition of ${\mathcal{L}}_{t}$, where $t\in \gf_p$, and
orthogonality relations, we have
$$pD_a=\sum_{t=0}^{p-1}{\mathcal{L}}_{t}\omega_p^{-at},$$
for any $a\in\gf_p$. It follows that in $\bC[(\gf_{p^n},+)]$, we
have
\begin{eqnarray}
pD_{\RR}&=&\bar{r}_0{\mathcal{L}}_{\RR}+\bar{n}_0{\mathcal{L}}_{\NN}+\frac{p-1}{2}\gf_{p^n},\\
pD_{\NN}&=&\bar{n}_0{\mathcal{L}}_{\RR}+\bar{r}_0{\mathcal{L}}_{\NN}+\frac{p-1}{2}\gf_{p^n}.
\end{eqnarray}

We now divide the proof into two cases.

\noindent{\bf Case 1.} $p\equiv 1$ (mod 4). In this case, $\RR$ is a
$(p,\frac {p-1}{2}, \frac{p-5}{4}, \frac{p-1}{4})$ PDS in the
additive group of $\gf_p$, we have the following equalities in
$\bZ[(\gf_{p},+)]$:
\begin{eqnarray}
\RR^2&=&\frac{p-5}{4}\RR+\frac{p-1}{4}\NN+\frac{p-1}{2},\\
\NN^2&=&\frac{p-1}{4}\RR+\frac{p-5}{4}\NN+\frac{p-1}{2},\\
\RR\NN&=&\frac{p-1}{4}(\gf_p\setminus\{0\}).
\end{eqnarray}
We now compute ${\mathcal L}_{\RR}^2$. From definition we have
${\mathcal L}_{\RR}^2=\sum_{t,s\in\RR}{\mathcal L}_t{\mathcal
L}_s=\sum_{t,s\in\RR, t+s\neq 0}{\mathcal L}_t{\mathcal
L}_s+\sum_{t\in\RR}{\mathcal L}_t{\mathcal L}_{-t}.$ Now using Part
(1) and (2) of Lemma~\ref{ProdL}, and noting that $n$ is assumed to
be even, we have
$${\mathcal L}_{\RR}^2=u\sqrt{p}^n\sum_{t,s\in\RR, t+s\neq 0}{\mathcal
L}_v+\frac{p-1}{2}p^n,$$ where
$v=(s^{1-\ell}+t^{1-\ell})^{\frac{1}{1-\ell}}$. Note that as $s$
(respectively, $t$) runs through $\RR$, so does $s^{1-\ell}$
(respectively, $t^{1-\ell}$). Now using (2.5), we have
\begin{equation}\label{lr2}
{\mathcal L}_{\RR}^2=u\sqrt{p}^n\bigg(\frac{p-5}{4}{\mathcal
L}_{\RR}+\frac{p-1}{4}{\mathcal L}_{\NN}\bigg)+\frac{p-1}{2}p^n.
\end{equation}
Similarly, we have
\begin{eqnarray}
{\mathcal L}_{\NN}^2&=&u\sqrt{p}^n\bigg(\frac{p-1}{4}{\mathcal
L}_{\RR}+\frac{p-5}{4}{\mathcal L}_{\NN}\bigg)+\frac{p-1}{2}p^n,\\
{\mathcal L}_{\RR}{\mathcal
L}_{\NN}&=&u\sqrt{p}^n\bigg(\frac{p-1}{4}{\mathcal
L}_{\RR}+\frac{p-1}{4}{\mathcal L}_{\NN}\bigg).
\end{eqnarray}
Since $p\equiv 1$ (mod 4), we have $\bar{r}_0=r_0$ and
$\bar{n}_0=n_0$. Now (2.3) becomes $pD_{\RR}=r_0{\mathcal
L}_{\RR}+n_0{\mathcal L}_{\NN}+\frac{p-1}{2}\gf_{p^n}.$ It follows
that
\begin{eqnarray*}
(pD_{\RR}-\frac{p-1}{2}\gf_{p^n})^2&=&(r_0{\mathcal
L}_{\RR}+n_0{\mathcal L}_{\NN})^2\\
&=&r_0^2{\mathcal L}_{\RR}^2+n_0^2{\mathcal
L}_{\NN}^2+2r_0n_0{\mathcal L}_{\RR}{\mathcal L}_{\NN}
\end{eqnarray*}
Using (\ref{lr2}), (2.9) and (2.10), we have
$$(pD_{\RR}-\frac{p-1}{2}\gf_{p^n})^2=\frac{p^2-1}{4}p^n+u\sqrt{p}^n(pD_{\RR}-\frac{p-1}{2}\gf_{p^n}).$$
By (\ref{eq1.3}), we see that $D_{\RR}$ is a PDS in $(\gf_{p^n},+)$.
Similar computations show that $D_{\NN}$ also satisfies the same
equation in $\bZ[(\gf_{p^n},+)]$, i.e.,
$$(pD_{\NN}-\frac{p-1}{2}\gf_{p^n})^2=\frac{p^2-1}{4}p^n+u\sqrt{p}^n(pD_{\NN}-\frac{p-1}{2}\gf_{p^n}).$$
Therefore $D_{\NN}$ is also a PDS in $(\gf_{p^n},+)$.

Now note that $pD_0={\mathcal L}_0+{\mathcal L}_{\RR}+{\mathcal
L}_{\NN}.$ So $$p^2D_0^2={\mathcal L}_0^2+{\mathcal
L}_{\RR}^2+{\mathcal L}_{\NN}^2+2{\mathcal L}_0{\mathcal
L}_{\RR}+2{\mathcal L}_0{\mathcal L}_{\NN}+2{\mathcal
L}_{\RR}{\mathcal L}_{\NN}.$$ Using (\ref{lr2}), (2.9), (2.10) and
the fact that ${\mathcal L}_0=\gf_{p^n}$, we have
$$p^2D_0^2=(-p^n+2p|D_0|-u\sqrt{p}^n(p-2))\gf_{p^n}+up(p-2)\sqrt{p}^nD_0+(p-1)p^n.$$
This shows that $D_0$ is a PDS in $(\gf_{p^n},+)$. Hence
$D_0\setminus\{0\}$ is also a PDS in $(\gf_{p^n},+)$.\\

\noindent{\bf Case 2.} $p\equiv 3$ (mod 4). In this case, $\RR$ is a
$(p,\frac{p-1}{2},\frac{p-3}{4})$ difference set in $(\gf_{p},+)$.
Hence we have the following equalities in $\bZ[(\gf_{p},+)]$:
\begin{eqnarray*}
\RR^2&=&\frac{p-3}{4}\RR+\frac{p+1}{4}\NN,\\
\NN^2&=&\frac{p+1}{4}\RR+\frac{p-3}{4}\NN,\\
\RR\NN&=&\frac{p+1}{4}+\frac{p-3}{4}\gf_p.
\end{eqnarray*}
By computations similar to those in Case 1, we now have
\begin{eqnarray*}
{\mathcal L}_{\RR}^2&=&u\sqrt{-p}^n\bigg(\frac{p-3}{4}{\mathcal
L}_{\RR}+\frac{p+1}{4}{\mathcal L}_{\NN}\bigg),\\
{\mathcal L}_{\NN}^2&=&u\sqrt{-p}^n\bigg(\frac{p+1}{4}{\mathcal
L}_{\RR}+\frac{p-3}{4}{\mathcal L}_{\NN}\bigg),\\
{\mathcal L}_{\RR}{\mathcal
L}_{\NN}&=&u\sqrt{-p}^n\bigg(\frac{p-3}{4}{\mathcal
L}_{\RR}+\frac{p-3}{4}{\mathcal L}_{\NN}\bigg)+\frac{p-1}{2}p^n.
\end{eqnarray*}
Now we can compute $(pD_{\RR}-\frac{p-1}{2}\gf_{p^n})^2$ by using
(2.3) and the above three equalities. We have
$$(pD_{\RR}-\frac{p-1}{2}\gf_{p^n})^2=\frac{p^2-1}{4}p^n+u\sqrt{-p}^n(pD_{\RR}-\frac{p-1}{2}\gf_{p^n}).$$
Similarly, we have
$$(pD_{\NN}-\frac{p-1}{2}\gf_{p^n})^2=\frac{p^2-1}{4}p^n+u\sqrt{-p}^n(pD_{\NN}-\frac{p-1}{2}\gf_{p^n}).$$
and
$$p^2D_0^2=(-p^n+2p|D_0|-u\sqrt{-p}^n(p-2))\gf_{p^n}+up(p-2)\sqrt{-p}^nD_0+(p-1)p^n.$$
We have shown that $D_0$, $D_{\RR}$ and $D_{\NN}$ are all PDS in
$(\gf_{p^n},+)$. The proof is now complete. \qed

Most of the known examples of $p$-ary weakly regular bent functions
are either quadratic forms or ternary (i.e., $p=3$). One can find a
table in \cite{TPF09} which summarizes the known examples at that
time. After \cite{TPF09} appeared, Helleseth and Kholosha
\cite{khHelleseth} constructed a class of $p$-ary weakly regular
bent functions $f$ from $\gf_{p^{4m}}$ to $\gf_p$, where $p$ is an
arbitrary odd prime and $f$ is not a quadratic form on
$\gf_{p^{4m}}$. One can apply Theorem~\ref{SC1} to this new class of
$p$-ary weakly regular bent functions to obtain partial difference
sets. Here we only want to demonstrate the applicability of
Theorem~\ref{SC1} and we are not concerned with the inequivalence
issues of the PDS constructed.

Next we prove that the PDS in Theorem~\ref{SC1} lead to amorphic
association schemes. The quickest way to prove this result is to
invoke the following theorem from \cite{vanDam, vdMuzy}.

\begin{theorem}{\em (\cite{vanDam, vdMuzy})} \label{3srg}
Let $\{G_1,G_2,G_3\}$ form a strongly regular decomposition of the
complete graph. Then $\{G_1,G_2,G_3\}$ forms an amorphic 3-class
association scheme.
\end{theorem}

Combining Theorem~\ref{SC1} and Theorem~\ref{3srg}, the following
corollary is immediate.
\begin{corollary}
Let $f:\gf_{p^{n}}\rightarrow \gf_{p}$ be a weakly regular bent
function. Assume that the Walsh coefficients of $f$ satisfy
(\ref{WRf}) and there exists a constant $k$ with $\gcd(k-1,p-1)=1$
such that $f(tx)=t^{k}f(x)$ for all $t\in\gf_{p}$ and all $x\in
\gf_{p^{n}}$. Furthermore, assume that $n$ is even. Then the three
Cayley graphs ${\rm Cay}(G,D_{0}\setminus \{0\})$, ${\rm
Cay}(G,D_{\RR})$ and ${\rm Cay}(G,D_{\NN})$, where
$G=(\gf_{p^n},+)$, form a 3-class amorphic association scheme on
$\gf_{p^{n}}$.
\end{corollary}

\vspace{0.2in}
\noindent{\bf Acknowledgement.} After we finished this paper in Jan. 2010, we were informed that 
Chee, Tan and Zhang \cite{chtz} also discovered the construction of PDS from weakly regular $p$-ary bent functions in Section 3 independently.

\end{document}